\documentclass[11pt]{article}
\usepackage{amsmath}
\usepackage{amsmath,amssymb,amsthm,latexsym,amscd,
indentfirst}

\title{\textbf{The  Space-like Surfaces with Vanishing Conformal
 Form
in the Conformal Space}
\footnote{Faculty of Mathematics and Computer Sciences,
 Hubei University,
 Wuhan 430062,
 People's Republic of China
 (e-mail: chxnie@163.com).}
\thanks {This work was partially supported by National Natural Science Foundation of China
  (Grant Nos.   10971055, 10801006) and
  Zhongdian Natural Science Foundation of Hubei Educational Committee. }
}
\author {   Changxiong  Nie }

\date{ }

\begin{document}
\maketitle
\noindent
{\bf Abstract.}
The conformal geometry of surfaces in the conformal space $\mathbf Q^n_1$ is studied. We classify the space-like
surfaces in $\mathbf Q^n_1$ with vanishing conformal form up to conformal equivalence.


\par\noindent
{\bf  1 Introduction.  }
\par\medskip

 In [1] Wang gave the structure
equations for M\"obius geometry of submanifolds in the unit sphere. Three fundamental
 tensors $\mathbb A$, $\mathbb B$ and $\Phi$ arise naturally in the structure equations.
 In [1]   $\mathbb A$ is called the Blaschke tensor,  $\mathbb B$  the M¡§obius second fundamental form
   , and   $\Phi$ the M\"obius form.
   Together
with M\"obius metric $g$, these tensors determine the submanifold up to M\"obius transformations of the unit sphere.

   Li and Wang$^{[2]}$ classified   surfaces with vanishing
 M\"obius form in sphere space $\mathbb S^{n+1}$. Readers should be reminded that Bryant$^{[3]}$ have
 classified all minimal surfaces with constant   curvature in
  the unit sphere $\mathbf S^n$, the hyperbolic space $\mathbf H^n$
 and Euclidean space $\mathbf R^n$. Li and Wang used Bryant's results.
 It is interesting to classify stationary surfaces with constant   curvature
in $\mathbf R^n_1,\mathbf S^n_1$ or $\mathbf H^n_1$.
Some other results about Lorentz conformal geometry see refs. [4-7]. Further relative knowledge
refers to   [8-10].

 We
conclude this paper with

 {\bf The Main Thoerem } Let $x:\mathbf M\rightarrow{\mathbf Q}^n_1 $ be
a regular spacelike full surface with vanishing form. Then $x$ are
the following four alternatives:

(i) $x$ is a stationary surface with constant   curvature
in $\mathbf R^n_1,\mathbf S^n_1$ or $\mathbf H^n_1$.

(ii) $x$ is a hyperbolic cylinder $\mathbf H^1\times\mathbf R$  in
$\mathbf R^3_1$.

(iii)  $x$ is a surface $\mathbf H^1(\sqrt{r})\times\mathbf
S^1(\sqrt{1+r})$ in $\mathbf S^3_1$, $r>0$.

(iv) $x$ is a hyperbolic torus $\mathbf H^1(\sqrt{1+r})\times\mathbf
H^1(\sqrt{-r})$ in $\mathbf H^3_1,-\frac{1}{2}\leq r< 0$.\\

\par\noindent
{\bf  2. The fundamental equations}
\par\medskip
 Let $\mathbb R^{n }_{ s}$ be the real vector space $\mathbb R^{n }$
with the Lorenzian inner product $\langle,\rangle$ given by
$$\langle X,Y\rangle:=\sum_{i= 1}^{n-s } x_iy_i
 - \sum_{i= n- s + 1 }^{n } x_iy_i,
   $$
  where $X=(x_1, \cdots
x_{n }),Y=(y_1, \cdots , y_{n })\in\mathbb R^{n }$.  We define the de
 Sitter sphere $\mathbf S^{n}_1$  and anti-de Sitter sphere
 $\mathbf H^{n}_1$ by
 $$\mathbf S^{n}_1=\{u\in \mathbf R^{n  +  1 }_1|\langle u, u\rangle=1\},
\ \ \ \mathbf H^{n}_1=\{u\in \mathbf R^{n +  1}_2|\langle u,
u\rangle=-1\}.   $$ We call Lorentzian space $\mathbf R^{n}_1$, de
Sitter sphere $\mathbf S^{n}_1$ and anti-de Sitter sphere
$\mathbf H^{n}_1$ Lorentzian space forms.
  Denote
   $${\mathbf Q}^n_1 =\{[x]\text{ is\ the\ projective\ coordinates}|
   \langle x,x\rangle=0,
  x\in \mathbf R^{n+2}_2\}.$$
 By some conformal diffeomorphisms we may regard    ${\mathbf Q}^{n}_1$ as the common compactified
space of $\mathbf R^{n}_1, \mathbf S^{n}_1$ and $ \mathbf
H^n_1$. In fact the conformal space ${\mathbf Q}^n_1 $ has a standard Lorentzian metric.
 We research the conformal geometry of surfaces under the conformal group of this
  Lorentzian metric. We refer the reader to [4, 7]   for further details.

 Suppose that $x:\mathbf M\rightarrow{\mathbf Q}^n_1 $ is a
space-like surface. That is, $x_*(\mathrm{T}\mathbf M)$
is non-degenerated subbundle of $\mathrm{T}{\mathbf Q}^n_1 $. Let
$y:U\rightarrow \mathbf R^{n+2}_2$ be a lift of $x:\mathbf
M\rightarrow{\mathbf Q}^n_1  $ defined in an open subset $U$ of
$\mathbf M$. We denote by $\Delta$ and $\kappa$    Laplacian   and
the normalized scalar curvature of the local positive definite
metric $\langle  d y,  d y\rangle $. Then we have

{\bf Theorem 2.1.} ( see [1] Theorem 1.2.)  On $\mathbf M$ the 2-form $g= - (\langle \Delta y,
\Delta y\rangle -4\kappa)\langle  d y,  d y\rangle   $ is a
globally defined invariant of $x:\mathbf M\rightarrow {\mathbf Q}^n_1   $ under the
Lorentz group transformations of ${\mathbf Q}^n_1  $.\\

Let $x:\mathbf M^2\rightarrow{\mathbf Q}^n_1 $ be a \textsl{regular} space-like
surface. That is,   the 2-form $g=-(\langle\Delta y,\Delta
y\rangle -4\kappa)\langle dy,dy\rangle $, which is called conformal metric,
 is  non-degenerated. Let $Y  =
 \sqrt{    -(\langle\Delta y,\Delta
y\rangle -4\kappa)      }  y$ be the canonical
lift of $x$ and define $N:\mathbf M\rightarrow\mathbf R^{n+2}_2$ by
$N= -\frac{1}{2}\Delta Y-\frac{1}{8}\langle\Delta Y,\Delta Y\rangle
Y$. Let $\{E_\alpha\}$ be a local basis of the conformal normal
bundle $\mathbf V$ of $x$. If $z=u+\mathrm{i}v$ be a local
isothermal coordinate on $\mathbf M$ for $g$, we can write
$g=e^{2\omega}|dz|^2=\frac{1}{2}e^{2\omega}(dz\otimes d\bar z+d\bar
z\otimes dz)$ for some local smooth function $\omega$. Denote by $K$
the Gauss curvature of $g$, we have $\Delta Y=4e^{-2\omega}Y_{z\bar
z},K=- 4e^{-2\omega}\omega_{z\bar z}$. Since
$\{Y,N,\text{Re}(Y_z),\text{Im}( Y_z),E_\alpha\}$ is a moving frame
in $\mathbf R^{n+2}_2$ along $\mathbf M$, one can
 write the structure equations and the fundamental equations as
 in [2]. Define
 $$\psi=2\langle N_z,Y_z\rangle ,\quad
 \phi_\alpha=\langle N_z,E_\alpha\rangle ,\quad
 \Omega_\alpha=2\langle Y_{zz},E_\alpha\rangle ,
 \quad A_{\alpha\beta}=\langle(E_\alpha)_z,E_\beta\rangle ,
 \eqno{(2.1)}$$
 and $ \phi^\alpha=\sum_{\beta}
g^{\alpha\beta}\phi_\beta,\quad\Omega^\alpha=\sum_{\beta}
g^{\alpha\beta}\Omega_\beta,\quad A^\beta_\alpha= \sum_\gamma
g^{\beta\gamma}A_{\alpha\gamma}$. The structure equations are
$$N_z=\frac{1}{8}(4K-1)Y_z+e^{-2\omega}\psi Y_{\bar z}+
\sum_\alpha \phi^\alpha E_\alpha,\eqno{(2.2)}$$
$$Y_{zz}=-\frac{1}{2}\psi Y+2\omega_z Y_{  z}+\frac{1}{2}
\sum_\alpha \Omega^\alpha E_\alpha,\quad
 Y_{z\bar z}=- \frac{1}{16}e^{2\omega}(4K-1)Y
-\frac{1}{2}e^{2\omega}N,\eqno{(2.3)}$$
$$(E_\alpha)_z=-\phi_\alpha Y-e^{-2\omega}\Omega_\alpha Y_{\bar z}+
\sum_\beta A^\beta_\alpha E_\beta.\eqno{(2.4)}$$ The fundamental
equations are
 $$\psi_{\bar z}=\frac{1}{2}e^{2\omega}K_z-
 \sum_\alpha\Omega^\alpha\bar\phi_\alpha, \quad
 \sum_\alpha\Omega^\alpha\Omega_\alpha=
 -\frac{1}{4}e^{4\omega}, \quad
 (\Omega_\alpha)_{\bar z}=-\sum_\beta
\Omega^\beta\bar
A_{\beta\alpha}-e^{2\omega}\phi_\alpha,\eqno{(2.5)}$$
$$(\phi_\alpha)_{\bar z}-\frac{1}{2}e^{-2\omega}
\bar\psi\Omega_\alpha+\sum_\beta\phi^\beta\bar A_{\beta\alpha} =
(\bar\phi_\alpha)_{ z}-\frac{1}{2}e^{-2\omega}
\psi\bar\Omega_\alpha+\sum_\beta\bar\phi^\beta A_{\beta\alpha},
\eqno{(2.6)}$$
$$(A_{\alpha\beta})_{\bar z}-(\bar A_{\alpha\beta})_{ z}=
\frac{1}{2}e^{-2\omega}(\Omega_\alpha\bar\Omega_\beta-\bar\Omega_\alpha\Omega_\beta)
+\sum_\gamma ( \bar
A_{\alpha\gamma}A^\gamma_\beta-A_{\alpha\gamma}\bar A^\gamma_\beta).
\eqno{(2.7)}$$

  {\bf Remark 2.1.} $\Psi=\psi
dz\otimes dz,\Phi=\sum_\alpha(\phi^\alpha dz + \bar \phi^\alpha
d\bar z)\otimes E_\alpha $ and $\Omega=\sum_\alpha\Omega^\alpha
dz\otimes dz\otimes E_\alpha$ are globally defined
conformal invariants.

 {\bf Remark 2.2.} The Willmore equations are
$$(\phi_\alpha)_{\bar z}-\frac{1}{2}e^{-2\omega}
\bar\psi\Omega_\alpha+\sum_\beta\phi^\beta\bar A_{\beta\alpha}
=0,\forall\alpha. \eqno{(2.8)}$$

\par\noindent
{\bf 3. The classification of space-like surfaces  in ${\mathbf Q}^n_1 $ with
$\Phi=0$}
\par\medskip
 Let $x:\mathbf M\rightarrow{\mathbf Q}^n_1 $ be a space-like surface in ${\mathbf Q}^n_1 $
 with vanishing conformal form, i.e.,
$\phi_\alpha=0,3\leq\alpha\leq n $. Then the fundamental equations
come into
$$\psi_{\bar z}=\frac{1}{2}e^{2\omega}K_z, \quad
\bar\psi\Omega_\alpha=\psi\bar\Omega_\alpha ,\quad
e^{-4\omega}\sum_\alpha\Omega^\alpha\Omega_\alpha=-
\frac{1}{4},\eqno{(3.1 )}$$
$$(\Omega_\alpha)_{\bar z}=-\sum_\beta
\Omega^\beta\bar A_{\beta\alpha},\quad
A_{\alpha\beta}=-A_{\beta\alpha}. \eqno{(3.2)}$$ It follows from
(3.2) that
$$(\sum_\alpha\Omega^\alpha\Omega_\alpha)_{\bar z}=0, \eqno{(3.3)}$$
thus the globally defined 4-form $\sum_\alpha
\Omega^\alpha\Omega_\alpha dz^4$ is holomorphic on $\mathbf M$. From
(3.1) we get
$$ \bar\psi\sum_\alpha\Omega^\alpha\Omega_\alpha
=\sum_\alpha\Omega^\alpha\bar\Omega_\alpha=-\frac{1}{4}
e^{4\omega}\psi.\eqno{(3.4)}$$
 Immediately we have

 {\bf Lemma 3.1.}  Let
 $x:\mathbf M\rightarrow{\mathbf Q}^n_1 $ be a surface in ${\mathbf Q}^n_1 $
 with vanishing conformal form. Then the conformal invariant
$\Psi=\psi dz^2$ is holomorphic on $\mathbf M$. Then $\Psi$ vanishes
identically or the zero points of  $\Psi$ are isolated.\\

First we consider the case that  $\Psi\equiv0$. Thus from (3.1) $K$
must be a constant. Then we get from (2.3) that
$$N=\frac{1}{8}(4K-1)Y+\mathbf c $$
for some constant vector $\mathbf c\neq0$ in $\mathbf R^{n+2}_2$.
Therefore $x$ is a conformal isotropic surface, i.e., $x$ is stationary
 surface with constant curvature in $\mathbf R^n_1,\mathbf S^n_1,$ or $\mathbf H^n_1$
 (see ref. [7]).

Now we come to discuss the case that the zero points of  $\Psi$ are
isolated. In this case we can cut $\mathbf M$ by some disjoint
curves $C_i$ to get a simply connected domain
 $\mathbf U=\mathbf M\backslash \sum_iC_i$
 such that $x:\mathbf U\rightarrow{\mathbf Q}^n_1 $ is a surface with $\Psi\neq0$ on $U$.
 By choosing complex
coordinate if necessary,
 we may assume that $\psi\equiv1$ on $\mathbf U$. It
 follows from (2.1) and (3.1) that $K=-4e^{-2\omega}\omega_{z\bar z}
 =0$. By (3.3) we know that $\{\Omega_\alpha\}$ are real functions.
 We define a global real vector field $E\in\mathbf V$ by
   $$ E=2e^{-2\omega}\sum_a\Omega^\alpha E_\alpha,\eqno{(3.5)}$$
then by (3.2) we have $\langle E,E\rangle =-1$. Choosing $\tilde
E_3=E$ and expanding it to a local orthonormal basis of the
conformal normal bundle $\{\tilde E_3,\tilde E_4, \cdots,\tilde
E_n\}$ , one can easily verify that
$$\tilde\Omega_3=-\frac{1}{2}e^{2\omega},\tilde\Omega_4=\cdots=\tilde\Omega
_{n}=0 .\eqno{(3.6)}$$ Using (3.3) , we get
$$(\tilde\Omega_3)_{\bar z}=-\sum^n_{\beta=3}
\tilde\Omega^\beta\bar A_{\beta3}=-\tilde\Omega_3\bar A_{33}=0,
\eqno{(3.7)}$$
$$0=(\tilde\Omega_\alpha)_{\bar z}=-\sum_{\beta\neq\alpha}
\tilde\Omega^\beta\bar A_{\beta\alpha}=-\tilde\Omega_3\bar
A_{3\alpha},\alpha> 3. \eqno{(3.8)}$$ Thus $\omega$ is a constant
and $A_{3\alpha}=0,\forall\alpha$. Now the structure equations read
$$N_z=\frac{1}{8}Y_z+e^{-2\omega}Y_{\bar z}, \eqno{(3.9)}$$
$$ Y_{zz}=-\frac{1}{2}Y-\frac{1}{4}e^{2\omega}\tilde E_3,\quad
 Y_{z\bar z}= \frac{1}{16}e^{2\omega}Y-\frac{1}
{2}e^{2\omega}N, \eqno{(3.10)}$$
$$(\tilde E_3)_z= \frac{1}{2}Y_{\bar z}, \quad
 (\tilde E_\alpha)_z=\sum_{\beta }A^\beta_\alpha \tilde E_\beta,\quad\alpha> 3. \eqno{(3.11)}$$
A surface is said to be full in ${\mathbf Q}^n_1 $ if $x(\mathbf M)$
does not lie in any totally umbilic $\mathbf Q^{n-1}$ of ${\mathbf Q}^n_1 $.
  We assume that $n\geq4$ . Then fixing a point $p\in U$, we
can find a constant vector $\xi\in\mathbf R^{n+2}_2$ with
$\langle\xi,\xi\rangle =1$ such that
$$\langle Y(p),\xi\rangle =0,\quad\langle N(p),\xi\rangle =0,\quad\langle Y_u(p),\xi\rangle
=0,$$
$$\langle Y_v(p),\xi\rangle =0,\quad\langle \tilde E_3(p),\xi\rangle =0. \eqno{(3.12)}$$ We
define real functions
$$f_1=\langle Y,\xi\rangle ,\quad f_2=\langle N,\xi\rangle ,\quad f_3=\langle Y_u,\xi\rangle ,$$
$$f_4=\langle Y_v ,\xi\rangle ,f_5=\langle \tilde E_3,\xi\rangle . \eqno{(3.13)}$$
Then  by (3.11)-(3.13) we can find constants $\{a_{\lambda\mu}\}$
and $\{b_{\lambda\mu}\}$ such that
$$(f_{\lambda})_u=\sum_\mu a_{\lambda\mu} f_\mu,\quad
(f_{\lambda})_v=\sum_\mu b_{\lambda\mu}
f_\mu,\quad1\leq\lambda,\mu\leq5. \eqno{(3.14)}$$ By (3.12) and the
uniqueness of the linear PDE (3.14) we get $f_\lambda \equiv0$. In
particular, $f_1=\langle Y,\xi\rangle =0$ on $U$, which implies that
$\langle Y,\xi\rangle =0$ on $\mathbf M$. If $x$ is full then $n$
must be 3.

\textsl{Acknownagements: The authors would like to express his
gratitude to Professor Changping   Wang
and Professor Haizhong  Li  for their warm-hearted instruction and inspiration. }\\

\noindent
 {\large\bf References}

\noindent
1. Wang C.P.,
  {Moebius geometry of submanifolds in $S^n$},
   Manuscripta Math.
{\bf 96}(1998): 517-534\\
2. Li H.Z.,   Wang C.P., {Surfaces with vanishing
Moebius form in $\mathbb{S}^n$}, Acta Math.  Sinica, Engl.
Ser., {\bf 19}(2003): 671-678\\
3. Bryant R.L.,
   {Minimal surfaces of constant curvature in
   $\mathbb S^n$},
   Trans. Amer. Math. Soc.
   {\bf 290}(1985): 259-271\\
4.   Deng Y.J., Wang C.P., Willmore surfaces in Lorentzian space, Sci. China Ser. A,
 {\bf 35}( 2005):
1361--1372\\
5.     Nie C.X., Ma X., Wang C.P.,
   Conformal CMC-surfaces in Lorentzian space forms,
 Chin. Ann. Math., Ser. B, {\bf 28} ( 2007): 299--310\\
6.   Alias L.J., Palmer  B., Conformal geometry of surfaces in Lorentzian space forms, Geometriae
Dedicata, {\bf 60}( 1996): 301--315\\
7.   Nie C.X., Li T.Z., He Y.J., Wu C.X., {Conformal isoparametric
  hypersurfaces with two distinct conformal principal
  curvatures in conformal space}, Sci. China Ser. A, {\bf 53} (2010): 953-965\\
8.   Blaschke W., Vorlesungen \"uber Differentialgeometrie, Vol. 3, Springer, Berlin, 1929\\
9.   Hertrich-Jeromin U., Introduction to M\"obius Differential Geometry, London Math. Soc. Lecture Note
Series, Vol. 300, Cambridge University Press, Cambridge, 2003\\
10.   O'Neill B., Semi-Riemannian Geometry with Applications to Relativity, Pure and Applied Mathematics,
103, Academic Press, New York, 1983

\end{document}